\documentclass[letterpaper, 10pt, conference]{ieeeconf}

\IEEEoverridecommandlockouts                            
\usepackage{amsmath,amsthm,enumerate,pdfsync,asymptote}  
\usepackage{setspace} 
\usepackage{graphicx,subfig}
\usepackage{amssymb}
\usepackage{setspace}
\usepackage{graphics,framed}
\usepackage{pgf}
\usepackage{tikz}
\usepackage{pgfplots,pgfmath}

\newtheorem{Definition}{Definition}
\newtheorem{Example}{Example}

\newtheorem{Remark}{Remark}
\newcommand{\R}{\mathbb{R}}

\newcommand{\Cp}{\mathbb{C}\mathrm{P}}

\usetikzlibrary{calc,fadings,decorations.pathreplacing}

\newcommand{\scone}[3]{%
\begin{scope}[rotate=#3,xshift=#1,yshift=#2]
\def\mypath{ (.-.2,0) -- +(-.2,.8) arc (180:0:.8) -- +(-.2,-.8) arc (-180:0:.2) } 
\fill [gray] mypath;
}

\tikzset{%
  >=latex, 
  inner sep=0pt,%
  outer sep=2pt,%
  mark coordinate/.style={inner sep=0pt,outer sep=0pt,minimum size=3pt,
    fill=black,circle}%
}

\makeatletter
\newsavebox{\sfe@box}
{\color@endgroup\egroup\subfloat[\sfe@caption]%
{\usebox{\sfe@box}}}
\makeatother

\title{\LARGE \bf
Known unknowns, unknown unknowns and information flow: new concepts and challenges in decentralized control}

\author{M.-A. Belabbas 
\thanks{M.-A. Belabbas is with the School of Engineering and Applied Sciences, Harvard
University, Cambridge, MA 02138 {\tt\small belabbas@seas.harvard.edu}}%
}

\begin{document}

\maketitle

\begin{abstract}             
We introduce and analyze a model for decentralized control. The model is broad enough to include problems such as formation control, decentralization of the power grid and flocking. The objective of this paper is twofold. First, we show how the issue of decentralization goes beyond having agents  know only part of the state of the system. In fact, we argue that  a complete theory of decentralization should take into account the fact that agents  can be made aware of only part of the global objective of the ensemble. A second contribution of this paper is the introduction of a rigorous definition of information flow for a decentralized system: we show how to attach to a general nonlinear decentralized system a unique \emph{information flow graph} that is an  invariant of the system.  In order to address some finer issues in decentralized system, such as the existence of so-called "information loops", we further refine the information flow graph to a simplicial complex---more precisely, a Whitney complex.   We illustrate the main results on a variety of examples.

\end{abstract}

\section{Introduction}

Informally speaking, a decentralized control system is a system whose different parts---let us call the different parts agents---are not told what to do by a unique, centralized controller, but decide what to do based on the possibly \emph{incomplete information} that is at their disposal. 

The importance of decentralization in control has been recognized for many decades~\cite{safonov78}, but only in more recent time has the issue been the subject of sustained investigation, see~\cite{lall06, shah09, KrishnamurthyKS08, murrayfax04, yu09,olfati06} and references therein. This renewed interest is fuelled, on the one hand, by the potential a complete theory of decentralized control has in explaining natural behavior: flocks of birds, ant colonies, etc. and more broadly by its applications in decision theory~\cite{SZjaamas08} and cognition. On the other hand, a theory of decentralization is also a necessity for  engineering design: from smart grids to  vehicles management on the highway~\cite{varaiya94,BarooahMehtaHespanhaAug09},  recent developments in robotics and communication have made it possible to envision very large groups of autonomous vehicles or agents collaborating to achieve a global objective. In this context,  decentralization  is thus necessary  for reasons ranging from robustness--- failure at some level (agent, controller, etc.) in a centralized system is likely to affect the entire system, whereas failure  in a decentralized system is more easily handled---to, at a more fundamental level, feasibility. Indeed,  a centralized controller for, say,  vehicles on the highway is not easily  implemented.

The objective of this paper is twofold. First, we will show that the idea of \emph{incomplete information} in a decentralized setting should be explored  beyond the usual partial knowledge of the state of the system. While most of the extant work in decentralized control implicitly assume that all the agents know  the objective of the ensemble and are thus solely constrained by their limited observations, we develop here a model which includes restrictions on what  agents know about the global objective of the ensemble. We illustrate this idea on a formation control problem in~\cite{belabbas_icca2011_mathform} by showing that local stabilization around a given configuration is possible only if agents know more than their selfish objective.

Second, we will make  rigorous the notion of \emph{information flow in a decentralized system}. The naive notion of information flow that is often used in the linear theory of decentralized system---i.e. splitting variables into groups  and coding the dependence between groups by a graph, see Section~\ref{sec:infflow}---is inherently dependent on the choice of coordinates used. This aspect puts it at odds with the idea that what an agent  knows about the system should not depend on the way one  chooses to describe the system.  Even more, the naive information flow does not acknowledge the possibility that Lie brackets may be needed to make the system controllable.  While these issues can be sidestepped in the linear case to obtain results that are nevertheless meaningful, they become  a genuine limitation when one tries to understand nonlinear decentralized systems or systems with constraints. Indeed, in the former situation, there often does not exist  preferred coordinates  or one may need several coordinate charts to describe the system. In the latter situation, choosing coordinates that are compatible with the constraints, e.g. a conservation of energy constraint, changes the naive information flow since it often requires a mixing of the coordinates used to describe the agents individually.  

We introduce in this paper a definition of  information flow graph  that is invariant under changes of coordinates and allows to define rigorously decentralization in a control system. The main idea behind this definition is that the observation functions on the system provide a natural set of vertices for the information flow graph.

The paper is organized as follows. We first introduce a general nonlinear model for decentralized control systems. We then define the \emph{global objective} of a system and the  \emph{local or selfish objectives} of the agents. The fact that the agents only know part of the global objective  is enforced through the use of non-invertible functions on the parameters describing the global objective; this approach can be understood, informally speaking, as a \emph{decentralization of the objective} or a decentralization of the design of the control law.

In the following section, we introduce a partial order on the local objectives and observation functions. This partial order allows us to quantify the idea that some local objectives (resp. observations) are more revealing of the global objective (resp. state of the ensemble) than others.

In Section~\ref{sec:infflow}, we introduce a coordinate free definition of information flow in a decentralized system. Motivated by the existence of "information loops" in  decentralized system,  we further refine the information flow graph into an information flow complex that reveals  finer issues in decentralization.  

\begin{Example}[Power grid]
 In recent years, the development of methods to insure the stability of the power grid have come at the forefront of research in control theory. In this context, one can view the power grid as a very large scale  system whose global objective is to remain stable around a desired operating point. Such systems currently operate in a centralized manner under the supervision of  Supervisory Control and Data Acquisition system (SCADA). This centralized framework, however,  is starting to show its limitations due to the increasingly complex components that are part of the grid (e.g. green energy suppliers). 

The development of decentralized methods to insure the stability and good operation of the grid have thus become a priority in power systems engineering. In this context, the global objective is a function of all the components of the grid, but it is clearly not feasible to let all agents in the grid know about its complete architecture.
\end{Example}

\begin{figure}
\begin{center}
\subfloat[Naive information flow for the $x$ variables]{
\begin{tikzpicture}[scale=.9] 

\node [fill=black,circle, inner sep=1pt,label=90:$x_1$] (1) at (1,1) {};
\node [fill=black,circle, inner sep=1pt,label=180:$x_2$] (2) at (0,0) {};
\node [fill=black,circle, inner sep=1pt,label=-90:$x_3$] (3) at (1,-1) {};
\node [fill=black,circle, inner sep=1pt,label=0:$x_4$] (4) at (2,0) {};

\draw [-stealth, ] (1) -- (2) node [ midway,above,  sloped, blue] {}; 
\draw [-stealth, ] (4) -- (3)  node [ midway,above,  sloped, blue] {};
\draw [-stealth, ] (2) -- (3)  node [ midway,above,  sloped, blue] {};
\draw [-stealth, ] (3) -- (1) node [ midway,above,  sloped, blue] {}; 
\draw [-stealth, ] (1) -- (4)  node [ midway,below,  sloped, blue] {};

\end{tikzpicture}\label{fig:infflow2cx}}\qquad
\subfloat[Naive information flow for the $z$ variables]{
\begin{tikzpicture}[scale=1.5] 

\node [fill=black,circle, inner sep=1pt,label=-135:$z_1$] (1) at (0,0) {};
\node [fill=black,circle, inner sep=1pt,label=180:$z_2$] (2) at (0,1) {};
\node [fill=black,circle, inner sep=1pt,label=90:$z_3$] (3) at (.5,1.5) {};
\node [fill=black,circle, inner sep=1pt,label=0:$z_4$] (4) at (1,1) {};
\node [fill=black,circle, inner sep=1pt,label=-45:$z_5$] (5) at (1,0) {};

\draw [-stealth, ] (1) -- (2) node [ midway,above,  sloped, blue] {}; 
\draw [-stealth, ] (1) -- (5) node [ midway,above,  sloped, blue] {}; 

\draw [-stealth, ] (2) -- (3) node [ midway,above,  sloped, blue] {}; 

\draw [-stealth, ] (3) -- (1)  node [ midway,above,  sloped, blue] {};
\draw [-stealth, ] (3) -- (5)  node [ midway,above,  sloped, blue] {};
\draw [-stealth, ] (4) -- (3) node [ midway,above,  sloped, blue] {}; 
\draw [-stealth, ] (5) -- (1)  node [ midway,below,  sloped, blue] {};
\draw [-stealth, ] (5) -- (4)  node [ midway,above,  sloped, blue] {};

\end{tikzpicture}\label{fig:infflow2cz}}

\end{center}\caption{The naive information flow of a system is a directed graph whose vertices $v_i$ correspond to  groups of variables describing the system.  There is a directed edge from $v_i$ to $v_j$ if the dynamics of variables in the group of $v_i$ depend on the variables in the group of $v_j$. This graph depends on the coordinates chosen to describe the system. }\label{fig:5agentsform}
\vspace{-.0cm}
\end{figure}
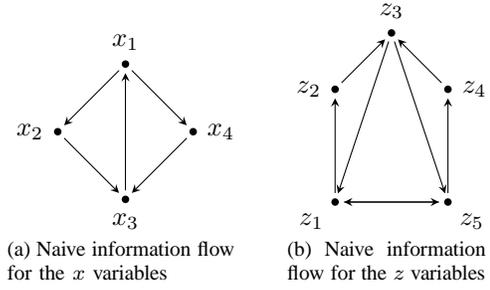

\begin{Example}[Formation control] 

Let $x_i \in \R^2$ represent the positions of autonomous agents in the plane and $d_i \in (0, \infty)$ be real positive constants. 

We consider the formation control problem whose dynamics are given by
\begin{eqnarray*}
\dot x_1 &=& e_1 (x_2-x_1)+e_5 (x_4-x_1)\\
\dot x_2 &=& e_2 (x_3-x_2)\\
\dot x_3 &=& e_3 (x_1-x_3)\\
\dot x_4 &=& e_4 (x_3-x_4)
\end{eqnarray*}
where we denote by $e_i$ the error in edge length:
\begin{multline*}
e_1 = \|x_{2}-x_1\|^2-d_1, e_2 = \|x_{3}-x_2\|^2-d_2, \ldots, \\
e_4 = \|x_{3}-x_4\|^2-d_4, e_5 = \|x_{4}-x_1\|^2-d_5 \end{multline*}

The  information flow of the system is represented in Figure~\ref{fig:infflow2cx}. Formation control problems are defined up to a rigid transformation of the plane~\cite{belabbascdc11sub1}. For this reason, one often describes the dynamics in terms of  the inter-agent distances~\cite{krick08}
\begin{equation}
\label{eq:defz}
\left\lbrace \begin{array}{rcl}
z_1 &=& x_2-x_1 \\
z_2 &=& x_3-x_2 \\
z_3 &=& x_1-x_3\\
z_4 &=& x_3-x_4 \\
z_5 &=& x_4-x_1,
\end{array}\right.
\end{equation}
instead of the absolute positions of the agents.

In the $z$ variables,  the dynamics of the system are given by 
\begin{eqnarray*}
\dot z_1 &=&  e_2 z_2-  e_1 z_1 -e_5 z_5\\
\dot z_2 &=& e_3z_3- e_2 z_2\\
\dot z_3 &=&  e_1 z_1+e_5z_5- e_3z_3\\
\dot z_4 &=&  e_3 z_3- e_4z_4\\
\dot z_5 &=&  e_4z_4-e_1 z_1-e_5z_5
\end{eqnarray*}

The corresponding naive information flow has 5 vertices and is represented in Figure~\ref{fig:infflow2cz}. We describe in Section~\ref{sec:infflow} a way to obtain the information flow depicted in Figure~\ref{fig:infflow2cx} from the description of the system given in terms of the $z$ variables.

\end{Example}

\section{A model for decentralized control}\label{sec:genmodel}

We present in this section a general model for nonlinear decentralized control systems. The model is a natural extension of  the notion of decentralized system that is encountered in the literature on linear systems.

In addition to being applicable to nonlinear problems, such as formation control, our approach distinguishes itself from most of the work on linear decentralized control in at least two major aspects~\cite{lall06,shah09,murrayfax04}:

\begin{itemize}
\item[-] it introduces the notion of parametrized objective of a decentralized system.
\item[-] it allows for \emph{loops of information} in the system, unlike approaches based on quadratic invariance or partial orders~\cite{lall06, shah09}. We revisit this point in Section~\ref{sec:infflow}.
\end{itemize}

\subsection{General model}

Let $M$ be a smooth manifold and the state $x \in M$. We consider nonlinear control systems of the type \begin{equation}\label{eq:defdecent1}\dot x =f(x,u(x))= \sum_{i=1}^n u_i(\delta_i(\mu); h_i(x)) g_i(x) \end{equation} where $\delta, h$ are smooth functions, the $g_i$'s are smooth vector fields and $\mu$ is a parameter that describes the objective of the system. We let $\mathcal U$ be the space of admissible controls $u_i$. We elaborate on the various parts of the model in this section.

A common situation is for the manifold $M$  to be the product of the manifolds describing the state-spaces of each agent: $$M = \bigotimes_{i=1}^n M_i$$ where $M_i$ is the state-space of agent $i$. We can thus write that the tangent space of $M$ is the direct sum $TM = \oplus_i TM_i$. In this  case, we also have that the projection of $g_i(x)$ onto $TM_j$ is zero if $i \neq j$: $$\pi_j g_i(x) = 0\mbox{ if }i \neq j.$$ 
This product structure is often lost due to either interactions between agents which impose constraints on the state $x\in M$,	  or the existence of a symmetry group acting on  $M$, in which case one has to consider equivalence classes of states $x \in M$ (this is the case in, e.g., formation control). Hence we do not assume here any special structure for $M$. 	

We differentiate between two type of objectives:
\begin{enumerate}
\item the objective that each agent or plant tries to satisfy: it is referred to as local objective or selfish objective.

\item the objective the agents try to achieve by cooperating: it is referred to as  global objective or common objective.
\end{enumerate}

The functions $\delta_i$ in Equation~\eqref{eq:defdecent1} allow us to control how much an agent knows about the common objective of the ensemble;  in some sense, these functions introduce a partial observation on the objective of the ensemble, akin to the partial observation that the agents have on the state of  the ensemble. 
They  are described in more detail below; we start with the definition of local observations.

\subsection{Local observations}\label{ssec:locandnaive}

The main characteristic of  a decentralized control system is that the agents are only able to observe part of the state of the system. We introduce the functions \vspace{-.0cm}$$h_i(x): M \rightarrow \R^{k_i}, k_i\mbox{ a positive integer}\vspace{-.0cm}$$ to describe the observation of agent $i$ on the current state of the system.   We denote by $h_i(M)$ the image of $M$ under the map $h_i$.

\vspace{-.0cm}
 
\subsection{Local and global objectives}
\vspace{-.0cm}
 
We define in this section the  \emph{local and global objective of a decentralized system}. We consider the  case of objectives that depend on a parameter. This level of generality is often necessary to accurately model decentralized systems  whose dynamics are rich enough to accommodate parameter-varying objectives, such as  flocks of autonomous agents or power distribution systems.

\subsubsection{Global objective}

Let $P$ be a smooth manifold, we let $\mu \in P$ parametrize the global objective of the control system as follows:

\begin{Definition}[Global objective]\label{def:globobj} Given a decentralized control system $\dot x = f(x,u(x))$ of the type of Equation~\eqref{eq:defdecent1},  the \emph{global objective function}  is a differentiable function $$F(\mu;x,u): P \times M \times \mathcal U \rightarrow \R^d$$ with the convention that the objective is achieved  if the system is at $x^* \in M$ with $$F(\mu;x^*,u) = 0$$ for \emph{equality objectives} or $$F(\mu;x^*,u) \geq 0$$ for \emph{inequality objectives}, where the inequality is taken entry-wise.
\end{Definition}

The objective function can in general depend on $u$; this dependence is necessary if ones considers stabilization objectives. When a global objective is not parametric or does not depend on $u$ explicitly, we omit the dependence from the notation. We give a few  examples:

\begin{Example}[Rendez-vous] Consider a multi-agent system with two agents whose positions are given by $x_1 \in \R^m$ and $x_2 \in \R^m$. The global objective is to have the  agents meet. This objective does not depend on a parameter. We can encode it by \vspace{-.0cm}$$F(x_1,x_2) = -\|x_1-x_2\|^2.\vspace{-.0cm}$$ If we want the agents to reach a position such that they are at a given distance  $d$ from each other, we let $P = [0, \infty)$ and we use $$\vspace{-.0cm}F(d;x_1,x_2) = -(\|x_1-x_2\|^2 - d^2)^2.\vspace{-.0cm}$$
\end{Example}

\begin{Example}[Stabilization] Consider the simple nonparametric rendez-vous problem described in the previous example with the addition that the agents are required to stabilize at the rendez-vous configuration. We denote by $$\frac{\partial f}{\partial x}|_{x^*}$$ the Jacobian  of the system at $x^*$.  We denote by $\lambda_i(A)$ the eigenvalue of $A$ with $i^{\mbox{th}}$ largest  real part.  We can represent this global objective by using the vector-valued function \begin{multline}
F(d;x,u) : \R^m \times \R^m\times \mathcal U\rightarrow \R^{m+1}: \\x \rightarrow \left[\begin{array}{c} -(\|x_1-x_2\|^2 - d^2)^2\\ -\mbox{Re}(\lambda_1(\frac{\partial f}{\partial x})) \\ \vdots \\ -\mbox{Re}(\lambda_n(\frac{\partial f}{\partial x}) ) \end{array}\right]\end{multline}
\end{Example}

\subsubsection{Local Objectives}

We now focus on describing the system at the level of the  agents. We define a \emph{local objective} as being, roughly speaking, a restriction of a global objective. 

\begin{Definition}[Local objective] Given a decentralized control system with global objective parametrized by $P$, we let 
$P_i$ be a smooth manifold and $$\delta_i : P \rightarrow P_i$$ be smooth functions. For $\mu \in P$,  the local objective of agent $i$ is given by a smooth function $$F_i(\delta_i(\mu);h_i(x),u_i): P_i\times h_i(M) \times \mathcal U_i \rightarrow \R^d$$ with the same convention as in Definition~\ref{def:globobj} regarding equality and inequality objectives.
\end{Definition}

When $\delta_i$ is not an invertible function, an agent \emph{knows about part of the global objective }. We give some examples of relations between global and local objectives in the section below.

The decentralized control problem is \emph{well-posed} if satisfying the local objectives is sufficient to satisfy the global objective: \begin{equation*}f_i(\delta_i(\mu),h_i(x),u_i(x)) \geq 0 \mbox{ for all } i  \implies F(\mu,x,f(x)) \geq 0\end{equation*} and similarly for the equality objective.

A wide array of questions in decentralized control can then be reduced to one on the following three major questions:

\begin{enumerate}
\item How little information can we let the agents know about the global objective and still have the ensemble achieve it? In other words, how informative do we need the $\delta_i$ to be in order to achieve a given global objective?

\item Given a global objective, how little observation on the system do the agents need in order to achieve the global objective? In other words, how informative do the $h_i$ need to be in order to achieve a global objective?

\item Given a decentralized system with fixed observation functions $h_i(x)$ and control vector fields $g_i(x)$, what global objectives   are achievable?
\end{enumerate}

The first two questions are  not  independent. Indeed, if the observation functions $h_i(x)$ do not provide much information about the state of the ensemble,  increasing the knowledge an agent has about the global objective is likely to be fruitless (a typical example is formation control). We introduce below a partial order on observations and objective that allow us to attach a mathematically precise meaning to these questions.

\begin{figure}
\begin{center}
\subfloat[]{
\begin{tikzpicture}[scale=.7]

\begin{scope}
\path[clip] (2-.5,1) -- (5-.5,1)--(6-.5,3) -- (5-.5,5) -- (2-.5,5) -- (1-.5,3) --cycle;

 \foreach \x in {1,2,...,5}
\foreach \y in {2,4} {
\node [] (\x\y) at (\x, \y) {};
} ;

\foreach \x in {1,2,3,...,5}
\foreach \y in {1,5} {
\node [] (\x\y) at (\x-.5, \y) {};
} ;

\foreach \x in {1,2,3,...,6}
\foreach \y in {3} {
\node [] (\x\y) at (\x-.5, \y) {};
} ;

\foreach \x in {1,2,3,...,5}
\foreach \y in {1,2,...,5} 
{
\draw [ very thin,densely  dotted  ] (\x\y) -- +(1,0);
\draw [ very thin,densely dotted  ] (\x\y) -- +(.5,1);
\draw [ very thin,densely dotted  ] (\x\y) -- +(.5,-1);
} 
\end{scope}

 \foreach \x in {1,2,...,5}
\foreach \y in {2,4} {
\node [fill=black,circle, inner sep=1pt,label=225:{\tiny $x_{\x\y}$}] (\x\y) at (\x, \y) {};
} ;

\foreach \x in {2,3,...,5}
\foreach \y in {1,5} {
\node [fill=black,circle, inner sep=1pt,label=225:{\tiny $x_{\x\y}$}] (\x\y) at (\x-.5, \y) {};
} ;

\foreach \x in {1,2,3,...,6}
\foreach \y in {3} {
\node [fill=black,circle, inner sep=1pt,label=225:{\tiny $x_{\x\y}$}] (\x\y) at (\x-.5, \y) {};
} ;
\end{tikzpicture}}
\subfloat[]{
\begin{tikzpicture}[scale=.7]

\begin{scope}
\path[clip] (2-.5,1) -- (5-.5,1)--(6-.5,3) -- (5-.5,5) -- (2-.5,5) -- (1-.5,3) --cycle;

 \foreach \x in {1,2,...,5}
\foreach \y in {2,4} {
\node [] (\x\y) at (\x, \y) {};
} ;

\foreach \x in {1,2,3,...,5}
\foreach \y in {1,5} {
\node [] (\x\y) at (\x-.5, \y) {};
} ;

\foreach \x in {1,2,3,...,6}
\foreach \y in {3} {
\node [] (\x\y) at (\x-.5, \y) {};
} ;

\foreach \x in {1,2,3,...,5}
\foreach \y in {1,2,...,5} 
{
\draw [ very thin,densely dotted  ] (\x\y) -- +(.5,1);
\draw [ very thin,densely dotted  ] (\x\y) -- +(.5,-1);
} 
\end{scope}

 \foreach \x in {1,2,...,5}
\foreach \y in {2,4} {
\node [fill=black,circle, inner sep=1pt,label=225:{\tiny $x_{\x\y}$}] (\x\y) at (\x, \y) {};
} ;

\foreach \x in {2,3,...,5}
\foreach \y in {1,5} {
\node [fill=black,circle, inner sep=1pt,label=225:{\tiny $x_{\x\y}$}] (\x\y) at (\x-.5, \y) {};
} ;

\foreach \x in {1,2,3,...,6}
\foreach \y in {3} {
\node [fill=black,circle, inner sep=1pt,label=225:{\tiny $x_{\x\y}$}] (\x\y) at (\x-.5, \y) {};
} ;
\end{tikzpicture}}

\end{center}
\caption{Consider the formation control problem where  agents with positions $x_{ij}\in \R^2$ are required to stabilize at the configuration depicted above. We  let $\mu$ be the vector of all pairwise distances between agents at the desired configuration. We represent the function $\delta_i(\mu)$ by a graph with a vertex per agent, and an edge between the vertices $x_{ij}$ and $x_{kl}$ if  agents $x_{ij}$ and $x_{kl}$ know the  distance $\|x_{ij}-x_{kl}\|$ at the desired configuration.  In order to have the agents cooperate to stabilize at this configuration, we can let each agent know about the complete vector $\mu$ or only parts of it.  For example, letting  each agent know about the distance to its nearest neighbors, as shown above in $(a)$, allows them to reconstruct the desired configuration. In $(b)$, the agents are given less information about the global objective than in $(a)$, and one can see that respecting the pairwise distances depicted is not sufficient to reconstruct the desired configuration. Hence the $\delta_i$ in $(b)$ are not informative enough. }\label{fig:exform1}
\vspace{-.0cm}
\end{figure}
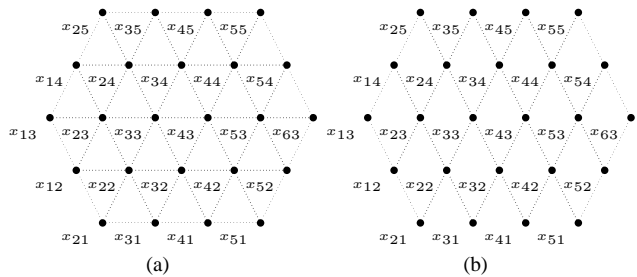

\begin{Example}[Formations]\label{ex:form1}
 Consider a formation control problem where  $n$  agents in the plane, with positions $x_i \in \R^2$, are required to stabilize at the configuration described in Figure~\ref{fig:exform1}. This configuration is defined up to a translation and rotation of the plane. We have shown in~\cite{belabbascdc11sub1} that the space of such configurations was $\Cp(n-2) \times (0,\infty)$. Hence, the parameter space is $P= \Cp(n-2)\times (0,\infty)$. A point in $P$ can also be represented, up to mirror symmetry, by the inter-agent distances $[d_1, \ldots, d_N]$, where $N=\frac{1}{2} n(n-1)$. The $d_i$'s are of course redundant in this representation, since simple trigonometric rules relate the pairwise distances. The problem of finding a non-redundant representation based on fewer pairwise distances is related to global rigidity~\cite{connelly05}

It may be impractical, or in some situation undesirable, to let every agent know about the complete vector $\mu$. The $\delta_i$ introduced here allow the study of systems where the amount of knowledge an agent gets about $\mu$ is controlled.

\end{Example}
\vspace{-.0cm}
\section{Partial orders and decentralization}
\vspace{-.0cm}
In the  design of a decentralized control system, a global objective can be achieved by the use of different observation functions and different local objectives.  We put a partial order on the local objectives (resp. observations) to formalize the notion that different local objectives (resp. observations)  can be more or less revealing of the global objective or system (resp. state of the ensemble).

We start with the definition of partial order:
\begin{Definition}[Partial order] A partial order $\succeq$ over a set $F$ is 	a binary relation on the elements of $F$ which satisfies, for $f_1,f_2,f_3 \in F$\vspace{-.0cm}
\begin{enumerate}
\item 	reflexivity: $f_1 \succeq f_1$.
\item antisymmetry: if $f_1 \succeq f_2$ and $f_2 \succeq f_1$ then $f_1 = f_2$
\item transitivity: $f_1 \succeq f_2$ and $f_2 \succeq f_3$ then $f_1 \succeq f_3$
\end{enumerate}\vspace{-.0cm}
\end{Definition}

The elements $f_1$ and $f_2$ in $F$ are called comparable if either $f_1 \succeq f_2$ or $f_2 \succeq f_1$ hold. The set $F$ is called a \emph{poset} for \emph{partially ordered set}. An element $f \in F$ is a greatest (resp. smallest) element if $f \succeq f_i$ (resp. $f_i \succeq f$) for all $f_i \in F$. If a greatest (resp. smallest) element exists, it is unique.

A maximal (resp. minimal) element $f$ is such that there is no $f_i \in F$ such that $f_i \succeq f$ (resp. $f \succeq f_i$).
\vspace{-.0cm}
\begin{Remark}
We mention here that partial orders have, quite interestingly, been applied to decentralized control in previous work~\cite{shah09}. However,  the approach and objective are quite different from ours. In the work~\cite{shah09}, the authors give an analysis of linear systems whose information flow graph---we will define it in Section~\ref{sec:infflow}---is given by a Hasse diagram~\cite{birkhoff_lattice48}. These are a type of \emph{directed acyclic} graphs.  They show in particular, relying on~\cite{lall06}, that there is a parametrization of such systems in which stabilization questions can be reduced to convex problems. 
\end{Remark}\vspace{-.0cm}
\vspace{-.0cm}
\subsection{Partial order on $\delta_i$}\label{ssec:pod}

\vspace{-.0cm}
Recall that the functions $\delta_i$ allow us to define decentralized systems where the  agents know only part of the global objective of the ensemble. We further refine this notion of incomplete knowledge of the global objective by establishing a \emph{partial order} of the functions $\delta_i$. Let $\mu \in P$, where $P$ is a smooth compact manifold. From Whitney's embedding theorem~\cite{warner83}, we know that for $n$ large enough, we can smoothly embed $P$ in $\R^n$. Hence, without loss of generality, we can assume that all the $\delta_i$ map into $\R^n$. We denote by $N_\delta(\mu)$ the isolevel set
$$N_\delta(\mu) = \lbrace x \in P \mbox{ s.t. } \delta(x) = \delta(\mu) \rbrace.$$

Many functions $\delta_i$ describe a similar restriction of the objective. For example, if $\delta_i$ maps to $\R$, translating the function by a constant $c \in \R$ to $\delta_i(\mu)+c$ does not change, for all practical purposes, what agent $i$ knows about the global objective. Indeed, if the objective is realized with the control $u_i$ for $\delta_i$, it is realized with the control $\tilde u_i(\delta_i;x)=u_i(\delta_i-c;x)$ for $\delta_i+c$. We generalize this idea in the following definition:
\vspace{-.0cm}
\begin{Definition}[Equivalence of local objectives]\label{def:equivrel} The functions $\delta_1:P \rightarrow \R^n$ and $\delta_2:P \rightarrow \R^n$ are equivalent at $\mu \in P$, written  $ \delta_1 \approx_\mu \delta_2$, if $$ N_{\delta_1}(\mu)=N_{\delta_2}(\mu).$$  They are equivalent if the above is true for all $\mu \in P$.
\end{Definition}\vspace{-.0cm}

Hence two functions are equivalent if their isolevel sets are the same. In particular, all one-to-one invertible functions are equivalent. This definition indeed corresponds to the intuitive notion of equivalent knowledge of the global objective:  as explained above, if the isolevel sets of $\delta_1$ and $\delta_2$  are the same, one can realize the same decentralized system by using an appropriately modified control $u$. 
\vspace{-.0cm}
\begin{Example}
Assume that we have a 2 agent decentralized system in $\R^m$ where the objective is parametrized by a point in $P=[0,1] \times [0,1]$. Let $\mu=(\mu_1,\mu_2) \in P$ and $\delta_1(\mu)=\mu_1$, $\delta_2(\mu)=\mu_2$. The uncertainty the first agent has about the global objective is its uncertainty about $\mu_2$. In particular, if agent 1 has access to $\tilde \delta_1(\mu)=\mu_1^2$, it has the same knowledge about the global objective than with $\delta_1(\mu)$.  According to Definition~\ref{def:equivrelf}, $\delta_1 \approx \tilde \delta_1$. Observe that this would not be true if $P=[-1,1]\times[-1,1]$.
\end{Example}\vspace{-.0cm}

We can now define a partial order on the local objectives.

\begin{Definition}[Partial order on $\delta_i$.]\label{def:parorderdelta} Let $F$ be the set of continuous functions on $P$ with the equivalence relation of Definition \ref{def:equivrel}. We say that $$\delta_1 \succeq_{\mu} \delta_2 \mbox{ if } N_{\delta_1}(\mu) \subseteq N_{\delta_2(\mu)}.$$  If the above relation is valid for all $\mu \in P$, we simply write $$\delta_1 \succeq \delta_2.$$

\end{Definition}

This definition expresses  the  notion that $\delta_1$ contains more information than $\delta_2$---written as $\delta_1 \succeq \delta_2$--- if the uncertainty arising from knowing $\delta_1(\mu)$ is smaller than the one arising from knowing  $\delta_2(\mu)$, where uncertainty is quantified by the isolevel sets of $\delta$. 

The partial order $\succeq$ has a smallest element: the constant function. This corresponds to the intuitive idea that if $\delta_i$ is constant, the  agent knows nothing about the global objective. The functions with the highest level of information are the invertible functions of $\mu$; these functions are the maximal elements. 
\vspace{-.0cm}

\begin{Example}Assume that $$P=S^3 = \lbrace x \in \R^4 \mbox{ s.t. } x_1^2+x_2^2+x_3^2+x_4^2 =1 \rbrace.$$
We let  $\delta_1(x)=x_1$ and $\delta_2(x)=(x_1,x_2)$.
We have that $$\delta_1 \succeq \delta_2.$$ We also have $$\delta_2 \approx (x_2,x_1) \approx (x_1+a, x_2+b), $$ for $a,b \in \R$
\end{Example}
\vspace{-.0cm}

\subsection{Partial order on $h_i$}
We similarly define a partial order on the observations $h_i$.  We denote by $N_{h_i}(x)$ the subset of $M$ such that $h_i(y) =h_i(x)$: $$N_{h}(x) = \lbrace y \in M \mbox{ s.t. } h(y)=h(x) \rbrace$$ In words, it is the set of  configurations  that are \emph{undistinguishable} to the observation function $h_i$.

Similarly to Definition~\ref{def:equivrel}, we say that two observation functions $h_1$ and $h_2$ are equivalent if their isolevel sets on $M$---or the configurations that are undistinguishable for $h_1$ and $h_2$---are the same:
$$h_1 \approx h_2 \Leftrightarrow N_{h_1}(x) = N_{h_2}(x),  \forall x \in M.$$

Furthermore, we can use a similar partial ordering on the observation functions to the one of Definition~\ref{def:parorderdelta}:$$h_1 \succeq_{x} h_2 \mbox{ if } N_{h_1}(x) \subseteq N_{h_2(x)}.$$  If the above relation is valid for all $x \in M$, we simply write $$h_1 \succeq h_2.$$

\vspace{-.0cm}

\subsection{Partial order on the  $f_i$}
We further define a partial order on the  $f_i$. In the case of equality objectives, the defintions are similar to the ones we have introduces above. We thus treat here the case of inequality objectives.  We denote by $N^+_{f_i(\delta_i;h_i,u_i)}(\mu)$ the subset of $M_i$ such that $f_i(\delta_i(\mu);h_i(x),u_i(x)) \geq 0$: $$N^+_{f_i(\delta_i;h_i,u_i)}(\mu) = \lbrace x \in M \mbox{ s.t. } f_i(\delta_i(\mu);h_i(x),u_i(x)) \geq 0 \rbrace$$ In words, it is the set of  configurations  that satisfy the local objective $f_i$.

We thus introduce the equivalence relation
\vspace{-.0cm}
\begin{Definition} 
\label{def:equivrelf} The functions $f_1$ and $f_2$ are equivalent at $\mu \in P$ if $$ f_1 \approx^+_\mu f_2 \Leftrightarrow N^+_{f_1(\delta_i,h_i,u_i)}(\mu)=N^+_{f_2(\delta_i,h_i,u_i)}(\mu).$$  They are equivalent if the above is true for all $\mu \in P$.

We say that $$ f_i \succeq_\mu f_j  \mbox { if } N^+_{f_i(\delta_i,h_i,u_i)}(\mu)\subseteq N^+_{f_j(\delta_i,h_i,u_i)}(\mu).$$
\end{Definition}
\vspace{-.0cm}
Hence, $f_1 \succeq f_2$ if the local objective $f_1$ is \emph{more stringent} than the local objective $f_2$.
\subsection{Minimally informed decentralized control and saturation}

We can now define
\vspace{-.0cm}
\begin{Definition}[Minimally informed decentralized system] Given a global objective $F$, we say that a decentralized control system of the type of Equation~\eqref{eq:defdecent1} is minimally informed if there is no set of functions $\tilde \delta_i, \tilde h_i, \tilde f_i$ with \vspace{-.1cm}$$ \delta_i \succeq \tilde \delta_i,   h_i \succeq \tilde h_i ,f_i \succeq \tilde f_i$$ and such that the decentralized system $$\dot x =\sum_i u_i(\tilde \delta_i(\mu),\tilde h_i(x)) g_i(x)\vspace{-.1cm}$$ with local objectives $\tilde f_i$ satisfies the global objective $F$.
\end{Definition}
\vspace{-.0cm}
 
When looking for a minimally informed system, an important notion that arises is the one of saturation.
Since the agents only have access to partial observations on the system, knowing an increasingly larger part of $\mu$ may cease to be helpful. 

We say that $\delta_i$ saturates $h_i$ if for all $\tilde \delta_i \succeq \delta_i$, there are no $\tilde f_i(\tilde \delta_i(\mu); h_i(x))$ with $$\tilde f_i \succeq f_i.$$ Reciprocally, we have that $h_i$ saturates $\delta_i$ if for all $\tilde h_i \succeq h_i$, there are no    $\tilde f_i( \delta_i(\mu); \tilde h_i(x))$ with $\tilde f_i \succeq f_i.$

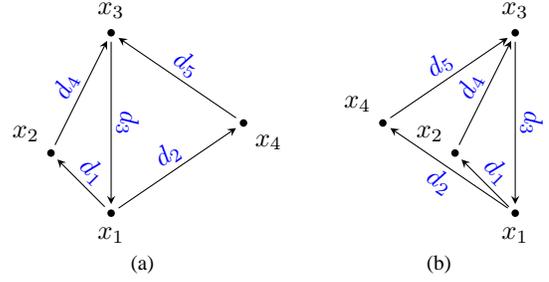
\begin{figure}
\begin{center}
\subfloat[]{
\begin{tikzpicture}[scale=.8,rotate=0] 

\node [fill=black,circle, inner sep=1pt,label=-90:$x_1$] (1) at ( 0, 0) {};
\node [fill=black,circle, inner sep=1pt,label=135:$x_2$] (2) at (-1 ,1) {};
  \node [fill=black,circle, inner sep=1pt,label=90:$x_3$] (3) at (0 , 3) {};
\node [fill=black,circle, inner sep=1pt,label=-45:$x_4$] (4) at ( 2.2,1.5) {};

\draw [-stealth,  ] (1) -- (2) node [ midway,above,  sloped, blue] {$d_1$}; ;
\draw [-stealth, ] (1) -- (4) node [ midway,above,  sloped, blue] {$d_2$};;
\draw [-stealth, ] (3) -- (1)  node [ midway,above,  sloped, blue] {$d_3$};;
\draw [-stealth, ] (2) -- (3)  node [ midway,above,  sloped, blue] {$d_4$};;
\draw [-stealth, ] (4) -- (3) node [ midway,above,  sloped, blue] {$d_5$};;
\end{tikzpicture}}\qquad
\subfloat[]{
\begin{tikzpicture}[scale=.8,rotate=0] 

\node [fill=black,circle, inner sep=1pt,label=-90:$x_1$] (1) at ( 0, 0) {};
\node [fill=black,circle, inner sep=1pt,label=135:$x_2$] (2) at (-1 ,1) {};
  \node [fill=black,circle, inner sep=1pt,label=90:$x_3$] (3) at (0 ,3) {};
\node [fill=black,circle, inner sep=1pt,label=135:$x_4$] (4) at ( -2.2,1.5) {};

\draw [-stealth,  ] (1) -- (2) node [ midway,above,  sloped, blue] {$d_1$}; ;
\draw [-stealth, ] (1) -- (4) node [ midway,below,  sloped, blue] {$d_2$};;
\draw [-stealth, ] (3) -- (1)  node [ midway,above,  sloped, blue] {$d_3$};;
\draw [-stealth, ] (2) -- (3)  node [ midway,above,  sloped, blue] {$d_4$};;
\draw [-stealth, ] (4) -- (3) node [ midway,above,  sloped, blue] {$d_5$};;
\end{tikzpicture}}

\end{center}
\caption{In the two-cycles formation depicted above, the agents are required to stabilize at the inter-agent distances $d_1, \ldots, d_5$. Up to mirror symmetry, there are two configurations in the plane that satisfy these interagent distances. If agent 1 can measure the \emph{relative positions} of agents 2 and 4, it can make use of the angle between the vectors $(x_2-x_1)$ and $(x_4-x_1)$ at the desired configurations $(a)$ and $(b)$. This observation function is thus not saturated by $d_1,d_5$. If agent 1 can only measure its distance to agent 1 and agent 2, the knowledge of the angle is not helpful. This observation function is saturated by $d_1,d_5$. } \label{fig:2cycle1}
\vspace{-.0cm}
\end{figure}

\begin{Example}\label{ex:2cycle1}
Consider agent 1 in the two-cycles formation of Figure~\ref{fig:2cycle1}. We let $\mu$ denote a target formation as explained in Example~\ref{ex:form1}. We let  $\delta_1(\mu) = [d_1, d_5]$. In the case of range only measurements  $$h_1(x) = \left[ \|x_2-x_1\|, \|x_4-x_1\|\right],$$ and the local objective is given by $$f_1(\delta_1(\mu);h_1(x)) = \left[\begin{matrix} \|x_2-x_1\|-d_1 \\ \|x_4-x_1\|-d_5 \end{matrix} \right].$$ It is easy to see that $h_1$ is saturated by $\delta_1$.  Similarly, for agent $2$ with $$h_2(x) = \left[ \|x_2-x_1\| \right],$$  $\delta_2(\mu) = [d_2]$ and $f_2(\delta_2(\mu);h_2(x)) =\|x_2-x_1\|-d_1$, we see that $\delta_2(\mu)=d_2$ saturates $h_2$.

If we let $$h_1(x) = \left[ \|x_2-x_1\|, \|x_4-x_1\|, (x_2-x_1)^T (x_4-x_1) \right],$$ i.e. the first agent observes its relative distances to agents 1 and 4 as well as their relative positions, then additional knowledge of $\mu$ is helpful. Indeed,  we can prove in this case that $h_1$ is saturated by $\delta_1(\mu) = \mu$. Intuitively, knowing all these distances allows agent $1$ to establish what the  possible angles between $x_2-x_1$ and $x_4-x_1$ are when  the global objective is reached. See~\cite{belabbas_icca2011_mathform} for additional details.

\end{Example}

\section{Information flow graph}\label{sec:infflow}

The abstract idea of information flow, because it allows to grasp the connectivity of the agents and identify potential difficulties in the distribution of information in the ensemble,  appears quite frequently in work on decentralized control.

We now define what we informally call the \emph{naive information flow} of a system; which is oftentimes implicitly defined in work on decentralized control. Let $e_j$ be the vector with zero entries except for the $j^{\mbox{th}}$ entry, which is one. The $e_j$'s form the canonical basis of  $\R^n$.  In the case of multi-agent systems in $\R^n$, $h_i(x)$ will often be the  projection of $x \in \R^n$ onto a the subspace spanned by some vectors $e_j$,  $j \in \mathcal J_i$ where $\mathcal J_i$ is a set of indices. For this reason, the observation functions $h_i$ are  encoded as a graph with vertices $x_i$ and an edge from $x_i$ to $x_l$ if $l \in \mathcal J_i$. We call this the \emph{naive} information flow of the system, since it is coordinate dependent as we illustrate in the examples below.

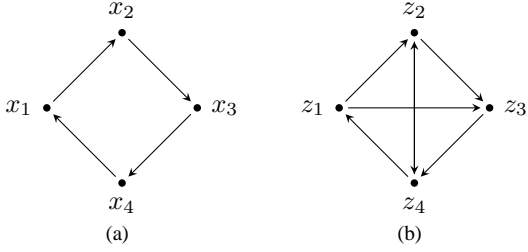
\begin{figure}[t]
\begin{center}
\subfloat[]{
\begin{tikzpicture} 

\node [fill=black,circle, inner sep=1pt,label=180:$x_1$] (1) at (0,0) {};
\node [fill=black,circle, inner sep=1pt,label=90:$x_2$] (2) at (1,1) {};
\node [fill=black,circle, inner sep=1pt,label=0:$x_3$] (3) at (2,.0) {};
\node [fill=black,circle, inner sep=1pt,label=-90:$x_4$] (4) at (1,-1) {};

\draw [-stealth,  ] (1) -- (2) node [ midway,above,  sloped, blue] {}; 
\draw [-stealth, ] (3) -- (4)  node [ midway,above,  sloped, blue] {};
\draw [-stealth, ] (2) -- (3)  node [ midway,above,  sloped, blue] {};
\draw [-stealth, ] (4) -- (1)  node [ midway,above,  sloped, blue] {};;

\end{tikzpicture}\label{fig:4agents1}}\qquad
\subfloat[]{
\begin{tikzpicture} 

\node [fill=black,circle, inner sep=1pt,label=180:$z_1$] (1) at (0,0) {};
\node [fill=black,circle, inner sep=1pt,label=90:$z_2$] (2) at (1,1) {};
\node [fill=black,circle, inner sep=1pt,label=0:$z_3$] (3) at (2,.0) {};
\node [fill=black,circle, inner sep=1pt,label=-90:$z_4$] (4) at (1,-1) {};

\draw [-stealth, ] (1) -- (2) node [ midway,above,  sloped, blue] {}; 
\draw [-stealth, ] (1) -- (3) node [ midway,above,  sloped, blue] {}; 

\draw [-stealth, ] (3) -- (4)  node [ midway,above,  sloped, blue] {};
\draw [-stealth, ] (2) -- (3)  node [ midway,above,  sloped, blue] {};
\draw [-stealth, ] (2) -- (4) node [ midway,above,  sloped, blue] {}; 

\draw [-stealth, ] (4) -- (1)  node [ midway,above,  sloped, blue] {};;
\draw [-stealth, ] (4) -- (2) node [ midway,above,  sloped, blue] {}; 

\end{tikzpicture}\label{fig:4agents2}}

\end{center}\caption{The two graphs above represent the naive information flow of the same system expressed in two different coordinate systems.}
\vspace{-.0cm}
\end{figure}

\vspace{-.0cm}
\begin{Example}[Four agents]\label{ex:4agents}
Consider the system with $M=\R^{4m}$, $x_i \in \R^m$ and whose dynamics is given by 
\vspace{-.0cm}
\begin{eqnarray*}
\dot x_1 &=& u_1(x_1,x_2) \\
\dot x_2 &=& u_2(x_2,x_3)\\
\dot x_3 &=& u_3(x_3,x_4)\\
\dot x_4 &=& u_4(x_4,x_1)
\vspace{-.0cm}
\end{eqnarray*}
Hence, $\mathcal J_1=\lbrace 1,2\rbrace, \mathcal J_2=\lbrace 2,3\rbrace,\mathcal J_3=\lbrace 3,4\rbrace,\mathcal J_4=\lbrace 1,4\rbrace.$ One can associate the graph of Figure~\ref{fig:4agents1} to this system, which shows a non-trivial loop in the information flow. Now consider the linear change of variables:
\vspace{-.0cm}
\begin{eqnarray*}
z_1 &=& x_1+x_2+x_3+x_4\\
z_2 &=& x_2+x_3+x_3\\
z_3 &=& x_3+x_4\\
z_4&=& x_4
\vspace{-.0cm}
\end{eqnarray*} We then have
\vspace{-.0cm}
\begin{eqnarray*}
x_1 &=& z_1-z_2\\
x_2 &=& z_2-z_3\\
x_3 &=& z_3-z_4\\
x_4 &=& z_4
\vspace{-.0cm}
\end{eqnarray*}
and for appropriately defined $\tilde u_i$,
\vspace{-.0cm}
\begin{eqnarray*}
\dot z_1 &=& \tilde u_1(z_1,z_2,z_3) \\
\dot z_2 &=& \tilde u_2(z_2,z_3,z_4)\\
\dot z_3 &=& \tilde u_3(z_3,z_4)\\
\dot z_4 &=& \tilde u_4(z_1,z_2,z_4).
\vspace{-.0cm}
\end{eqnarray*}

In this case, $\mathcal J_1=\lbrace 1,2,3\rbrace,\mathcal J_2=\lbrace 2,3,4\rbrace,\mathcal J_3=\lbrace 3,4\rbrace,\mathcal J_4=\lbrace 1,2,4\rbrace.$ This system corresponds  to the naive information flow graph depicted in Figure~\ref{fig:4agents2}. 
\end{Example}

\subsection{Information flow graph of a decentralized system}

We have seen  above that the naive definition of information flow is not satisfactory since it depends on the parametrization chosen for the system, whereas  decentralization is  a coordinate-free notion: our choice of coordinates to describe a system should not affect the knowledge each agent has about the system. 

We provide here a coordinate free definition of information flow. The idea is to let the observation functions $h_i$, define the vertices of the graph, and use the vector fields $g_i$ to determine the presence of edges,  Precisely, to a decentralized control system of the  type
\vspace{-.0cm}$$\dot x = \sum_{i}u_{i}(\delta_i(\mu);h_i(x))g_{i}(x)\vspace{-.0cm}$$ we will assign a directed graph at first, and refine the notion to obtain a \emph{simplicial complex}.

Recall that a vector field $g(x)$ on $M$ acts on functions $h$ defined on $M$ via differentiation. We write this action as $$g \cdot h(x).$$ If $h(x)=[h^1,\ldots, h^k]$ is vector-valued, we define $g \cdot h$ as $g \cdot h = [g \cdot h^1, \ldots g \cdot h^k]$. 

For example, on $\R^n$ with coordinates $(x_1,\ldots,x_n)$, the vector fields $g_1(x)=[g_{11}(x),\ldots,g_{1n}(x]$ and $g_2(x)=[g_{21}(x),\ldots, g_{2n}(x)]$ act on the function $h(x)$ via
$$g_i(x) \cdot h = \sum_{j=1}^n g_{ij} \frac{\partial}{\partial x_j} h .$$

The Lie bracket of $g_1$ and $g_2$ is the vector field $$[g_1,g_2](x) =\frac{\partial g_2}{\partial x} g_1 -\frac{\partial g_1}{\partial x} g_2  $$ where $\frac{\partial g}{\partial x}$ is the Jacobian matrix of $g$.

\begin{Definition}[Information Flow Graph]\label{eq:definfflow} Consider the decentralized control system \begin{equation}\label{eqdefinflow}\dot x = \sum_{i=1}^n \sum_{j=1}^{n_i} u_{ij}(\delta_i(\mu); h_i(x)) g_{ij}(x)\end{equation} where all the functions and vector fields involved are smooth. 
We assign to this system the graph with $n$ vertices $h_1, h_2, \ldots, h_n$ and edges given according to the following rules:
\begin{enumerate}
\item[$n_i=1$]
 
 there is an edge from $h_j$ to $h_i$ if $$g_{ik}(x) \cdot h_j(x) \neq 0 \mbox{ for any } k =1\ldots n_i$$ 
\item[$n_i \neq 1$] Let $\lbrace g_{i1},\ldots,g_{in_i}\rbrace_{LA}$ be the set of vector fields obtained by taking  iterated Lie brackets of $ g_{i1},\ldots,g_{in_j}$. There is an edge between from $h_j$ to $h_i$ if $$g_{ik}(x) \cdot h_j(x) \neq 0 \mbox{ for any } g_{ik}\in \lbrace g_{i1},\ldots,g_{in_i}\rbrace_{LA}.$$ 
\end{enumerate}

In words, there is an edge from $h_j$ to $h_i$ if the motion of an agent that uses the observation function $h_i$ is observable by $h_j$.
\end{Definition}

In the multi-agent case, each agent will often have its own observation function, and the above can be rephrased as saying that there is an edge from agent $i$ to agent $j$ if agent $i$ can observe changes in the state of agent $j$.

\begin{Example} Consider the system  on $\R^4$ given by
\begin{eqnarray*}\label{eq:3agexv}
\dot x_1 &=& u_1(x_1,x_2) \\
\dot x_2 &=& u_2(x_2,x_3)\\
\dot x_3 &=& u_3(x_3,x_4)\\
\dot x_4 &=& u_4(x_4,x_1)
\end{eqnarray*}

We let $x=(x_1,x_2,x_3,x_4)$ and  \begin{multline*}h_1(x)=(x_1,x_2), h_2(x)=(x_2,x_3), h_3(x)=(x_3,x_4)\\ \mbox{ and }h_4(x)=(x_1,x_4).\end{multline*}

We define \begin{multline*} g_1(x) = [1,0,0,0], g_2(x) = [0,1,0,0], \ldots,\\ g_4(x) =[0,0,0,1]\end{multline*} The system of Equation~\eqref{eq:3agexv} can thus be written as
$$\dot x = \sum_i u_i(h_i(x)) g_i(x).$$

According to definition~\ref{def:infflow}, the information flow graph is given by $G=(V,E)$ with $$V= \lbrace h_1, h_2,h_3, h_4\rbrace$$ and $$E=\lbrace (h_1,h_2), (h_2,h_3), (h_3,h_4), (h_4,h_1) \rbrace.$$

The same system expressed in the $z$ variables defined in Example~\ref{ex:4agents} is given by
\begin{eqnarray*}\label{eq:3agexv}
\dot z_1 &=& u_1(z_1-z_2,z_2-z_3)+u_2(z_2-z_3,z_3-z_4)\\ &&+u_3(z_3-z_4,z_4) +u_4(z_4,z_1-z_2)\\
\dot z_2 &=&u_2(z_2-z_3,z_3-z_4)+u_3(z_3-z_4,z_4)\\&&\hfill +u_4(z_4,z_1-z_2)\\
\dot z_3 &=&  u_3(z_3-z_4,z_4) +u_4(z_4,z_1-z_2)\\
\dot z_4&=&u_4(z_4,z_1-z_2).
\end{eqnarray*}

The observation functions are \begin{multline*}\tilde h_1(z)=(z_1-z_2,z_2-z_3), \tilde h_2(z)=(z_2-z_3,z_3-z_4), \\\tilde h_3(z)=(z_3-z_4,z_4), \tilde h_4(z)=(z_4,z_1-z_2).\end{multline*}The control vector fields become \begin{multline*}\tilde g_1(z)=[1,0,0,0], \tilde g_2(z)=[1,1,0,0], \tilde g_3(z)=[1,1,1,0]\\ \mbox{ and } \tilde g_4(z)=[1,1,1,1].\end{multline*}

We can now write $$\dot z =\sum_i u_i(\tilde h_i(z)) \tilde g_i(z).$$
The information flow graph associated to this system has four vertices $\tilde h_1, \tilde h_2, \tilde h_3, \tilde h_4$ We have the following relations
\begin{eqnarray*}
\tilde g_1 \cdot \tilde h_2 &=&[\partial_{z_1} (z_2-z_3), \partial_{z_1}(z_3-z_4) = [0,0] \\
\tilde g_1 \cdot \tilde h_3 &=&[0,0]; \mbox{    } \tilde g_1 \cdot \tilde h_4 =[0,1] \\
\tilde g_2 \cdot \tilde h_1 &=&[0,1]; \mbox{    }\tilde g_2\cdot \tilde h_3 =[0,0] \\
\tilde g_2 \cdot \tilde h_4 &=&[0,0]; \mbox{    }\tilde g_3 \cdot \tilde h_1 =[0,0] \\
\tilde g_3 \cdot \tilde h_2 &=&[0,1]; \mbox{    }\tilde g_3 \cdot \tilde h_4 =[0,0] \\
\tilde g_4 \cdot \tilde h_1 &=&[0,0]; \mbox{    }\tilde g_4 \cdot \tilde h_2 =[0,0] \\
\tilde g_4 \cdot \tilde h_3 &=&[0,1]
\end{eqnarray*}
These relations yield the same information flow graph as above.

\end{Example}

\subsection{Decentralized systems and Whitney complex}

Definition~\ref{def:infflow} in the previous section attaches  an information flow graph to a decentralized system in a coordinate free manner. The salient point was that the observation functions $h_i(x)$ provide a natural set of vertices for the graph.

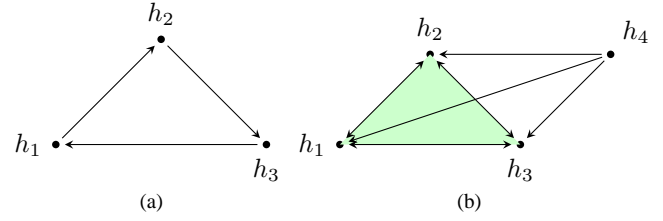
\begin{figure}
\begin{center}
\subfloat[]{\begin{tikzpicture}[scale=1.4]
\node [fill=black,circle, inner sep=1pt,label=180:$h_1$] (1) at ( 0, 0) {};
\node [fill=black,circle, inner sep=1pt,label=90:$h_2$] (2) at ( 1, 1) {};
\node [fill=black,circle, inner sep=1pt,label=270:$h_3$] (3) at ( 2, 0) {};

\draw [-stealth ] (1) -- (2);
\draw [-stealth ] (2) -- (3);
\draw [-stealth ] (3) -- (1);

\end{tikzpicture}\label{fig:triangform} }
\subfloat[]{
\begin{tikzpicture}[scale=1.2]
\node [fill=black,circle, inner sep=1pt,label=180:$h_1$] (1) at ( 0, 0) {};
\node [fill=black,circle, inner sep=1pt,label=90:$h_2$] (2) at ( 1, 1) {};
\node [fill=black,circle, inner sep=1pt,label=270:$h_3$] (3) at ( 2, 0) {};
\node [fill=black,circle, inner sep=1pt,label=45:$h_4$] (4) at ( 3, 1) {};

\fill [fill=green!20] (0,0)-- (1,1) -- (2,0)--cycle;
\draw [stealth-stealth ] (1) -- (2);
\draw [stealth-stealth ] (2) -- (3);
\draw [stealth-stealth ] (1) -- (3);
\draw [-stealth ] (4) -- (1);
\draw [-stealth ] (4) -- (2);
\draw [-stealth ] (4) -- (3);

\end{tikzpicture}\label{fig:illinfcomp}}
\end{center}\caption{In $(a)$, we represent the information flow complex of system~\eqref{eq:sysloop}, which exhibits a nontrivial information loop. In $(b)$, we represent the information flow complex of system~\eqref{eq:sysnoloop}; the shaded region depicts a 2-simplex in the complex. The information loop between $1,2$ and $3$ is trivial in this case.}
\end{figure}

Consider the triangular formation of Figure~\ref{fig:triangform}. The main source of difficulty in the control of this formation comes from the fact that the motion of $x_i$ depends on the motion of $x_{i+1}$ (taken modulo 3): if $x_2$ moves, $x_1$ has to adjust itself, which forces $x_3$ to move which in turn provokes a motion of $x_2$. We call this a nontrivial \emph{loop of information}.

Now consider the triangular formation with bidirectional edges. The above mentioned loop of information still exists, but its effect on the dynamics is diluted due to the fact that the communication goes both ways between the agents. In fact, there is all-to-all communication between agents in this formation, and the system is thus equivalent to a centralized one, where each agent implements locally a copy of a centralized controller. We call this information loop trivial. We devote the remainder of this section to putting this notion of triviality on a firm mathematical footing. 

In order to do so, we need some concepts from algebraic topology. The role of homological algebra and algebraic topology in control theory and applied sciences has been recognized in many different contexts such as feedback stabilization, computer graphics,  sensor networks or data analysis~\cite{brockett83asymptotic, Silva07homologicalsensor,edelsbrunner04,Ghrist07barcodes}. We show here how  related ideas naturally appear in the definition of decentralized systems.

We start with some graph theoretic definitions. We recall here that  the information flow graph is in general a mixed graph (i.e. containing both directed and undirected edges). We say that $G=(V,E)$ is an \emph{undirected complete graph} if $$E= \lbrace (v_i,v_j)\mbox{ s.t. } v_i, v_j \in V\rbrace.$$ In words, $G$ contains all possible undirected edges on its vertices. If $G=(V,E)$ is a graph, we call $G'=(V',E')$  the subgraph of $G$ generated by $V' \subset V$ when $E' \subset E$ is the set of edges of $E$ which start and end at  vertices in $V'$:
$$E' = \lbrace (v_i, v_j) \in E \mbox{ for all } v_i,v_j \in V'\rbrace.$$
A subgraph $G'$ of $G$ is an \emph{undirected clique} if it is  an undirected complete graph. We can now give a coordinate independent definition of decentralized systems:
\begin{Definition}[Decentralized system]
A system of type of Equation~\eqref{eq:definfflow} is \emph{centralized} if its associated information flow is an undirected complete graph. Otherwise, it is decentralized.
\end{Definition}

We now address the fact that some information loops are trivial, as described in the beginning of this section. A \emph{path} of length $k$ in a graph $G=(V,E)$ is an ordered list of vertices $v_1,\ldots, v_k$, without repetitions except possibly for $v_1$ and $v_k$, such that $(v_i,v_{i+1}) \in E$ for $i=1\ldots,k-1.$ A path is \emph{closed} or a \emph{loop} if $v_1=v_k$.

\begin{Definition}[Information loop]
A nontrivial information loop in a decentralized system with information flow graph $G=(V,E)$ is a closed path $(v_1,\ldots,v_k)$ such that the subgraph $G'$ generated by $V'=\lbrace v_1, \ldots, v_k \rbrace$ is not an undirected clique.
\end{Definition}

This definition takes into account the fact that when a graph is fully connected, even though loops will exist, their presence has no effect on the dynamics of the system.

The definition of information loop points towards the use of techniques from homological algebra to handle the information flow graph.  We define here a combinatorial object, called simplicial complex, which allows us to make a connection between the structure of decentralized systems and algebraic topology.

A k-simplex  is determined by k+1 vertices; we  use the usual notation  $[x_1,x_2,\ldots,x_{k+1}]$ for the k-simplex with vertices $x_1, \ldots, x_{k+1}$. A k-simplex has k+1 facets which are (k-1)-simplices, they are given by $[x_2,\ldots, x_k], [x_1,x_3,\ldots,x_k],\ldots,[x_1,\ldots,x_{k-1}]$. 
\begin{Definition}
An \emph{abstract simplicial complex} $\mathcal{S}$ is a combinatorial object consisting  of a set of simplices such that  any facet  of a simplex $s \in \mathcal S$ is also in $\mathcal S$. The \emph{k-skeleton} of a simplicial complex is the set of simplices of dimension $k$ or less. 
\end{Definition}

We have the following definition:
\begin{Definition}[Information Flow Complex]\label{def:infflow} Consider the decentralized control system $$\dot x = \sum_{i=1}^n \sum_{j=1}^{n_i} u_{ij}(\delta_i(\mu); h_i(x)) g_{ij}(x)$$ where all the functions and vector fields involved are smooth. We assign to this system the simplicial complex with $n$ vertices $h_1, h_2, \ldots, h_n$ and facets given according to the following:
\begin{enumerate}
\item There is an edge between $x_i$ and $x_j$ if $$g_{jk}(x) \cdot h_i(x) \neq 0 \mbox{ for any } k =1\ldots n_j$$ 
\item There is a k-simplex with vertices $x_1, \ldots, x_k$ if \begin{multline*} g_{jk} \cdot h_i \neq 0 \mbox{ for any } g_{jk} \in \lbrace g_{j1},\ldots,g_{jn_j} \rbrace_{LA}, \\ \mbox{ for all } i,j=1..k\end{multline*}
\end{enumerate}
\end{Definition}

The simplicial complex defined above is sometimes called a \emph{Whitney complex} or \emph{flag complex} in the literature. Due to space constraints, and the amount of background necessary to analyze such objects any further,  most notably via their cohomology groups, we leave the study of the information flow complex to future work.

\begin{Example}
Consider the system
\begin{equation}\label{eq:sysloop}
\left\lbrace \begin{matrix}
\dot x_1 &= &u_1(x_1,x_2)\\
\dot x_2 &= &u_2(x_2,x_3)\\
\dot x_3 &= &u_3(x_3,x_1)
\end{matrix}\right.
\end{equation}

Using the notation introduced above, we have
$$h_1(x)=(x_1,x_2), h_2(x)=(x_2,x_3), h_3(x)=(x_3,x_1)$$ and $g_i(x) = e_i.$
We thus have
\begin{eqnarray*}
g_1 \cdot h_2=[0,0] ;&& g_1 \cdot h_3=[0,1]\\
g_2 \cdot h_1=[0,1] ;&& g_2 \cdot h_3=[0,0]\\
g_3 \cdot h_1=[0,0] ;&& g_3 \cdot h_2=[0,1]
\end{eqnarray*}
We thus associate the information complex $$S=\lbrace h_1,h_2,h_3, [h_1,h_2], [h_2,h_3], [h_3,h_1]\rbrace.$$ 

Consider the system
\begin{equation}\label{eq:sysnoloop}
\left\lbrace\begin{array}{rcl}
\dot x_1 &= &u_1(x_1,x_2,x_3)\\
\dot x_2 &= &u_2(x_1,x_2,x_3)\\
\dot x_3 &= &u_3(x_1,x_2,x_3,x_4)\\
\dot x_4 &=& u_4(x_1,x_2,x_4)
\end{array}\right.
\end{equation}

Using the same approach as above, we find that the information complex of this system is (see  Figure~\ref{fig:illinfcomp})
 \begin{multline*}S=\lbrace h_1,h_2,h_3, h_4,[h_1,h_2], [h_2,h_3], [h_3,h_1], \\ [h_1,h_4],[h_2,h_4],[h_3,h_4], [h_1,h_2,h_3]\rbrace. \end{multline*}
\end{Example}

\bibliographystyle{IEEEtran}
\bibliography{distribcontrolbib}

\end{document}